\documentclass[12pt]{article}
\usepackage{latexsym,amsfonts,amsthm,amsmath,amscd,amssymb,color}

\usepackage{tikz}
\usepackage{tikz-cd}
\usetikzlibrary{matrix,arrows,decorations.pathmorphing}

\newtheorem{lemma}{Lemma}[section]

\newtheorem{rem}[lemma]{Remark}

\newtheorem{example}[lemma]{Example}
\newtheorem{defn}[lemma]{Definition}

\newcommand{\Ker}{\operatorname{Ker}}

\def\C{{\mathbb C}}
\def\R{{\mathbb R}}
\def\H{{\mathbb H}}

\def\L{{\cal L}}

\def\H{{\cal H}}
\def\G{{\cal G}}

\def\su{{\mathfrak{su}}}

\def\u{{\mathfrak{u}}}

\def\sl{{\mathfrak{sl}}}

\def\g{\mathfrak{g}}

\def\t{\tilde}
\def\t{\mathfrak{t}}

\def\u{{\mathfrak u}}
\def\su{{\mathfrak su}}

\def\L{{\cal L}}
\def\H{{\cal H}}
\def\t{{\bf t}}

\def\L{{\cal L}}

\def\R{{\mathbb R}}

\newcommand{\bea}{\begin{eqnarray}}
\newcommand{\eea}{\end{eqnarray}}

\newcommand{\Tr}{\textrm{Tr}}

\DeclareMathAlphabet{\mathpzc}{OT1}{pzc}{m}{it}

\newtheorem{theorem}{Theorem}[section]
\newtheorem{proposition}[theorem]{Proposition}
\newtheorem{corollary}[theorem]{Corollary}

\begin{document}

\title{Multiplicative integrable models \\ from Poisson-Nijenhuis structures}

\author{F. Bonechi \footnote{\small INFN Sezione di Firenze, email: francesco.bonechi@fi.infn.it}}


\maketitle

\begin{abstract}
We discuss the role of Poisson-Nijenhuis geometry in the definition of multiplicative integrable models on symplectic groupoids. These 
are integrable models that are compatible with the groupoid structure in such a way that the set of contour levels of the hamiltonians in involution
inherits a topological groupoid structure. We show that every maximal rank PN structure defines such a model. We consider the examples defined on compact hermitian 
symmetric spaces studied in \cite{BQT}.\footnote{Talk given at the conference ''{\it From Poisson Brackets to Universal Quantum Symmetries}'', August 18-22 2014, 
Stephan Banach International Mathematical Center, Warsaw}
\end{abstract}

\thispagestyle{empty}

\eject

\section{Introduction}

The notion of symplectic groupoid has its first motivation in the problem of quantizing Poisson manifolds. In fact, to any (integrable) Poisson manifold $M$ we can 
associate a canonical symplectic manifold $\G$ of double dimension endowed with a compatible groupoid structure. The compatibility is expressed by saying that the graph 
of the multiplication is a lagrangian submanifold of $\bar\G\times\bar\G\times\G$, where $\bar\G$ means that we consider the opposite symplectic structure. 
Any quantization scheme has to produce a space of states $\cal H$ and a state in $\cal H$ for any lagrangian submanifold. When quantization is applied to the graph of the 
multiplication of $\G$, we get a state $m\in\H^*\otimes\H^*\otimes\H$, that can be thought of as a map $m:\H\otimes\H\rightarrow\H$. Associativity of the groupoid 
multiplication implies that $m$ defines an algebra structure on $\H$, that must be considered as
the algebra of operators that quantizes the Poisson algebra of functions on $M$. So any quantization method of the symplectic groupoid should produce an algebra as output.

In \cite{BCST2,BCQT} we started a program to study the quantizatization of integrable Poisson manifolds through singular lagrangian polarizations of the symplectic groupoid.
The starting point, considered in \cite{Hawkins} first, is that in order to quantize a symplectic groupoid one has to look for polarizations that are compatible with the 
groupoid structure. Thanks to this compatibility, the application of geometric quantization methods should define a non commutative algebra on the 
space of polarized sections. The basic remark underlying \cite{BCQT} is that this perspective, rather than putting additional conditions on the difficult task of 
finding suitable polarizations, can allow more singular choices, namely singular real polarizations. Indeed, one can look for a singular real polarization such 
that the space of lagrangian leaves inherits the structure of topological groupoid. If this is the case the Bohr-Sommerfeld conditions select a subgroupoid that we call
as the {\it Bohr-Sommerfeld groupoid}; if it allows a Haar system, one can define the convolution algebra and regard this non commutative algebra as the output of the quantization
procedure. These conditions are completely different from those that make geometric quantization work in the real case, this requiring the space of leaves to be a smooth manifold.
In particular, one can consider as a source of possible examples integrable models on the symplectic groupoid. We gave then in \cite{BCQT} the definition of {\it multiplicative 
integrable model} as an integrable model whose space of contour levels inherits a groupoid structure. 

\medskip
\begin{defn}\label{def_multiplicative_integrable_model}
An integrable model $F\equiv \{f_i\}$ on the symplectic groupoid $\G$ is said to be multiplicative if the the contour level set $\L_F$ of the hamiltonians $F$ has a topological
groupoid structure such that the quotient $\G\rightarrow\L_F$ is a surjective groupoid morphism.
\end{defn}

The {\it modular class} is a distiguished class in Lichnerowicz-Poisson cohomology, that measures the existence of a volume form invariant with respect to hamiltonian
vector fields. A representative of this class depends on the choice of a volume form; such a representative is lifted to a groupoid one-cocycle of $\G$, that we call
the {\it modular function}. If the modular function is in involution with the hamiltonians $f_i$, then we say that the modular function is {\it multiplicatively integrable}.

This general framework was applied in \cite{BCQT} to the case of a family of Poisson structures on complex projective spaces. These are Poisson structures covariant with 
respect to the homogeneous action of the special unitary group, endowed with the standard Poisson-Lie structure. In this case, the integrable model comes from the canonical
hierarchy of a PN structure long ago introduced by \cite{KRR} for all compact hermitian symmetric spaces. In \cite{BCQT} we proved that the integrable model can be lifted
from the complex projective space to a multiplicative integrable model on the symplectic groupoid and computed the Bohr-Sommerfeld groupoid. 

The task of this note is to make clear the role of maximal rank Poisson Nijenhuis structures in the definition of multiplicative integrable models. We show that the 
construction of the multiplicative integrable model of \cite{BCQT} is actually valid in general for all maximal rank PN structures. 

In Section 2 we review the 
definition and the basic notions of PN structure. All the described results are well known; nevertheless, in the spirit of a concise but pedagogical introduction,
we included the proofs. 
This geometrical structure was introduced in \cite{MM} as a geometrical approach to integrable models. A PN structure is a Poisson structure $P$ compatible with a $(1,1)$-tensor $N$ with vanishing Nijenhuis torsion. As a consequence there is a full hierarchy
of compatible Poisson structures $P_{j+1}=N P_j $. The hamiltonians $I_k=\Tr N^k/k$ are in involution with respect to all $P_j$; we call them the {\it canonical hierarchy}.
If the first Poisson structure $P_1$ is non degenerate then the PN structure is said to be symplectic. The requirement that the hamiltonians of the canonical hierarchy are 
independent and so define
a completely integrable model leads to the definition of maximal rank PN structure. 
In particular we introduce the Poisson algebra of hamiltonian forms and show that it is abelian when the symplectic PN structure is of maximal rank. 

In Section 3 we prove the main result of this paper in Proposition \ref{pn_quotient_groupoid}: we show how the algebra of hamiltonian forms can be lifted to 
the symplectic groupoid integrating $P_2$ to a Poisson subalgebra; if the PN structure is of maximal rank
then this subalgebra is abelian and defines a multiplicative integrable model. This result explains the construction of \cite{BCQT}.

In \cite{KRR} it was shown that on compact hermitian symmetric spaces the inverse of the {\it Kirillov-Kostant-Souriau symplectic form} $\omega_{kks}$ and
the so-called {\it Bruhat-Poisson} structure $\pi$ are compatible and define a PN structure.
In \cite{BQT} we proved that these PN structures are of maximal rank. The review of this construction is the subject of Section 4. By applying our result, we conclude 
that there exists a multiplicative integrable model on the symplectic groupoid integrating $\pi+t\omega_{kks}^{-1}$ for each $t$. This model has been studied so far 
only for the case of complex projective spaces. Indeed it has been shown in \cite{BCQT} that in this case 
the multiplicative integrable model defines a lagrangian fibration and the Bohr-Sommerfeld groupoid has been computed.
The properties of the multiplicative integrable model on the symplectic groupoid in the general case will be the object of further investigation.

\bigskip
\bigskip

\section{From PN structures to integrable models}
In this section we recall basic facts about Poisson Nijenhuis geometry and its relation with integrable models. The results
are standard and can be found in \cite{MKS,MM,KS,DaFe,vaisman}. We reserve a particular attention to the notion of hamiltonian forms. For the sake
of completeness and to clarify certain points that are relevant to our construction, we give the explicit proofs of these results, that can be found
scattered in the literature.

A $(1,1)$-tensor $N:TM\rightarrow TM$ is called a {\it Nijenhuis} tensor if its Nijenhuis torsion $T(N)$ vanishes, 
{\it i.e.} for any
couple $(v_1,v_2)$ of vector fields on $M$ we have
\begin{equation}
\label{Njienhuis_torsion}
T(N)(v_1,v_2)=[N v_1,N v_2]-N ([N v_1,v_2] + [v_1,N v_2]-N [v_1,v_2] )=0  ~.
\end{equation}

It can be shown that the following anchor and 
bracket on $v,w \in {\rm Vect}^1(M)$ 
\begin{equation}\label{nijenhuis_algebroid_structure}
a_N(v) = N(v)\ ,\;\; [v,w]_N=[N(v),w] + [v,N(w)] - N([v,w]) \;\;\; \;.
\end{equation}
define an algebroid structure on $TM$ that we denote with $T_NM$. 
Remark that the vanishing torsion can be written as $0= [N(v_1),N(v_2)]-N([v_1,v_2]_N)$, so that we conclude that the bundle map $N:T_NM\rightarrow TM$ is an algebroid
morphism, where we denote with $TM$ the tangent algebroid. 
Let $\iota_N$ be the degree $0$ derivation on multivector fields defined as $\iota_N(f)=0$ and $\iota_N(X)=N(X)$ for $f\in C^\infty(M)$ and $X\in{\rm Vect}^1(M)$. Let 
$\iota_{N^*}$ be the dual derivation on $\Omega^1(M)$.
The algebroid differential is the degree one derivation $d_N$ on $\Omega^1(M)$ defined as
\begin{equation}
\label{Nijenhuis_algebroid_differential}
d_N = [\iota_{N^*},d]\;,
\end{equation}
that squares to zero. It is clear that $[d,d_N]=d d_N+d_N d = 0$. The {\it hamiltonian forms} are defined as 
\begin{equation}\label{definition_hamiltonian_forms}
\Omega^1_{ham}(M,N)=\{\alpha\in\Omega^1(M)|\, d\alpha=d_N\alpha=0\}\;.
\end{equation}
The name hamiltonian is justified when we look at them
in the context of the PN structures that we are going to introduce. Nevertheless their definition makes sense in this more general setting.

\smallskip
\begin{lemma}
The hamiltonian forms $\Omega^1_{ham}(M,N)$ are $N^*$-invariant.
\end{lemma}
{\it Proof}. Since $N:T_NM\rightarrow TM$ is an algebroid morphism then $N^*d=d_N N^*$ so that if $\alpha\in\Omega^1_{ham}(M,N)$ then $dN^*\alpha=d\iota_{N^*}\alpha=d_N\alpha=0$ 
and $d_N(N^*\alpha)=N^*d\alpha=0$ 
so that $N^*\alpha\in\Omega^1_{ham}(M,N)$. \qed

\smallskip
We recall that given a bivector $P\in\Gamma(\Lambda^2 TM)$ we can define the following antisymmetric bracket on $\Omega^1(M)$
\begin{equation}\label{bracket_one_forms}
\{\alpha,\beta\}_P = L_{P(\alpha)}(\beta) - L_{P(\beta)}(\alpha) - d\langle P,\alpha\wedge\beta\rangle\;\;\;  \alpha,\beta\in\Omega^1(M)\;, 
\end{equation}
where we denote with the same symbol the bivector $P$ and the antisymmetric map $P:T^*M\rightarrow TM$. This antisymmetric bracket satisfies the Jacobi identity and so is a Lie bracket if and only if $P$ is a Poisson bivector. In this case this bracket
and the anchor $P$ define an algebroid structure that we denote with $T^*_PM$. The algebroid differential is then computed as $d_P=[P,-]$, where the bracket is the Schouten bracket between multivectorfields, and its 
cohomology is called the Lichnerowicz-Poisson cohomology of $P$. Remark that $(TM,T^*_PM)$ is a bialgebroid, {\it i.e.} $d_P$ is a derivation of the Schouten bracket on 
multivector fields or, equivalently, the de Rham differential $d$ is a derivation of the bracket (\ref{bracket_one_forms}).

\smallskip
\begin{defn} 
A triple $(M,P,N)$, where $(M,P)$ is a Poisson manifold and $N$ a Nijenhuis tensor, is called a {\it Poisson-Nijenhuis} {\rm(PN)} manifold if
$P$ and $N$ are compatible, {\it i.e.}
\begin{equation}\label{compatibility_pnmfld}
N P = P N^*~,~~~~ \{\alpha,\beta\}_{NP}=\{N^*\alpha,\beta\}_P + \{\alpha,N^*\beta\}_P - N^*\{\alpha,\beta\}_P  ~~~~,
\end{equation}
for $\alpha,\beta\in \Omega^1(M)$.
\end{defn}

\medskip
The first of equations (\ref{compatibility_pnmfld}) says that $N P$ is antisymmetric and so defines a bivector; the bracket $\{,\}_{NP}$ appearing 
in the second of (\ref{compatibility_pnmfld}) is then the bracket between one forms defined by this bivector. The following characterization of PN manifolds has been 
proven in \cite{KS}.

\begin{proposition}
\label{PN_bialgebroid}
If $(M,P,N)$ is a {\rm PN} manifold then $(T_NM,T^*_PM)$ is a bialgebroid.
\end{proposition}
{\it Proof}. We prove first the identity
\begin{equation}\label{hierarchy_algebroid_differentials}
d_{NP} = [\iota_{N},d_P].
\end{equation}
Since (\ref{hierarchy_algebroid_differentials}) is an equality between derivations, it is enough to prove it on functions and vector fields. On functions it is a trivial check. 
Let $X\in{\rm Vect}^1(M)$ and $\alpha,\beta\in\Omega^1(M)$.
\begin{eqnarray*}
\langle \alpha\wedge\beta, d_{NP}(X)\rangle &=& NP(\alpha)\langle \beta,X\rangle - NP(\beta)\langle \alpha,X\rangle - \langle \{\alpha,\beta\}_{NP},X\rangle \cr
&=& PN^*(\alpha)\langle \beta,X\rangle - PN^*(\beta)\langle \alpha,X\rangle \cr
& & - \langle \{N^*\alpha,\beta\}_P + \{\alpha,N^*\beta\}_P-N^*\{\alpha,\beta\}_P,X\rangle\cr
&=& \langle N^*\alpha\wedge\beta,d_P(X)\rangle +P(\beta)\langle N^*\alpha,X\rangle \cr
& & + \langle \alpha\wedge N^*\beta,d_P(X)\rangle - P(\alpha)\langle N^*\beta,X\rangle + 
\langle\{\alpha,\beta\}_P ,NX\rangle \cr
&=& \langle \alpha\wedge\beta, [\iota_{N},d_P]X\rangle\;,
\end{eqnarray*}
where we used the definition of the algebroid differential in terms of the anchor and the Lie bracket and, in the second equality, properties (\ref{compatibility_pnmfld}).
Let $X,Y\in{\rm Vect}^1(M)$; we compute
\begin{eqnarray*}
d_P([X,Y]_N) &=& d_P([NX,Y] + [X,NY] -N[X,Y]) \cr
&=& [d_P(NX),Y] - [NX,d_P(Y)] + [d_PX,NY] - [X,d_P(NY)] \cr
& & -d_P(N[X,Y]) \cr
&=& [(\iota_{N}d_P-d_{NP})X,Y]- [NX,d_P(Y)]+ [d_PX,NY]\cr
& & -[X, (\iota_{N} d_P - d_{NP})(Y)] - (\iota_{N}d_P-d_{NP})([X,Y])\cr
&=& [\iota_{N}d_P(X),Y]+ [d_PX,NY]-\iota_{N}[d_PX,Y]\cr
& & - [NX,d_P(Y)]-[X,\iota_{N} d_P(Y)] + \iota_{N}[X,d_P(Y)] \cr
&=& [d_P(X),Y]_N - [X,d_P(Y)]_N\;,
\end{eqnarray*}
where we repeatedly used the fact that $d_P$ and $d_{NP}$ are derivations of the Schouten bracket between vector fields and, in the last equality,
we used a straightforward expression for the Gerstenhaber bracket between multivectorfields generated by the Lie bracket in (\ref{nijenhuis_algebroid_structure}) .
\qed

\smallskip

In the following Proposition we list the basic consequences of the above definition.

\smallskip
\begin{proposition}
\label{consequences_pn_structure} Let $(M,P,N)$ be a PN manifold.
\begin{itemize}
\item[$i$)] For all $r>0$ and $j\geq 0$, $(M,P_j,N^r)$ where $P_{j+1}= N^{j} P$, is a $PN$ manifold. Moreover, they are compatible, 
{\it i.e.} $[P_j,P_s]=0$.
 \item[$ii$)] The bundle map $P_k$ is an algebroid morphism $P_k:T^*_{P_{r}}M\rightarrow T_{N^{r-k}}M$ for each $r>k$.
 \item[$iii$)] The hamiltonian forms $\Omega^1_{ham}(M,N)$ define a $N^*$-invariant subalgebra with respect to the Koszul bracket $\{,\}_{P_j}$ of all the Poisson structures $P_j$ 
of the PN hierarchy. It is valid the following property for all $\alpha,\beta\in\Omega^1_{ham}(M,N)$
\begin{equation}
\label{exchange_rule}
N^*\{\alpha,\beta\}_{P_j}=\{N^*\alpha,\beta\}_{P_j} = \{\alpha,N^*\beta\}_{P_j}\;, 
\end{equation}
\end{itemize}
\end{proposition}
{\it Proof}.
($i$) Theorem 1.3 of \cite{vaisman}.

($ii$) The bundle map $P$ clearly intertwines the anchors of $T^*_{NP}M$ and $T_NM$. Let us show that it is a Lie algebra morphism. Indeed,
\begin{eqnarray*}
P(\{\alpha,\beta\}_{NP}) &=& P(\{N^*\alpha,\beta\}_P + \{\alpha,N^*\beta\}_P-\{N^*\alpha,N^*\beta\}_P) \cr
&=& [PN^*(\alpha),P(\beta)] + [P(\alpha),PN^*(\beta)]-[PN^*(\alpha),PN^*(\beta)] \cr
&=& [NP(\alpha),P(\beta)] + [P(\alpha),NP(\beta)]-[NP(\alpha),NP(\beta)] \cr
&=& [P(\alpha),P(\beta)]_N\;,
\end{eqnarray*}
where in the first and third equality we used (\ref{compatibility_pnmfld}) and in the second one we used the fact that $P$ is an algebroid morphism between $T^*_PM$ and $TM$.
By applying it to the PN structure $(M,P_k,N^{r-k})$ we get the general result.

($iii$) Let $\alpha,\beta\in\Omega^1_{ham}(M,N)$. From the definition (\ref{bracket_one_forms}), $\{\alpha,\beta\}_P$ is an exact form.
Since $(T_NM,T^*_PM)$ is a bialgebroid then $d_N$ is a derivation of the bracket defined by $P$ so that $d_N\{\alpha,\beta\}_P=0$. The result for $\{\alpha,\beta\}_{P_j}$ follows
by applying this result to the PN structure $(M,P_j,N)$

We compute, by using the first of (\ref{compatibility_pnmfld}), 
\begin{eqnarray*}
\{\alpha,\beta\}_{NP} &=& d\langle NP,\alpha\wedge\beta\rangle = d\langle NP(\alpha),\beta\rangle = d \langle P(\alpha),N^*(\beta)\rangle = \{\alpha,N^*(\beta)\}_P \cr
&=& d\langle PN^*, \alpha\wedge \beta\rangle = \langle PN^*(\alpha),\beta\rangle= \{N^*(\alpha),\beta\}_P\;.
\end{eqnarray*}
By applying the second of (\ref{compatibility_pnmfld}), we get (\ref{exchange_rule}). The same property for $P_j$ is easily obtained by applying it to $(M,P_j,N)$. \qed

\medskip

We denote as $\g_j(M,N)$ the Lie algebra structure on $\Omega^1_{ham}(M,N)$ defined by $P_j$. Since $(T_NM,T^*_{P_j}M)$ is a Lie bialgebroid, $d_N$
is a derivation of the Gerstenhaber bracket on forms defined by $P_j$ and, in particular, $H^1(T_NM)$ inherits the Lie algebra structure from each $P_j$ and the image 
of $\g_j(M,N)$ is a Lie subalgebra.


The following proposition clarifies the connection between $PN$ structures and integrable models.

\smallskip
\begin{proposition}\label{PN_hierarchies}
Let $\alpha\in\Omega^1_{ham}(M,N)$ and let $\alpha_{j+1}=N^*{}^j(\alpha)\in\Omega^1_{ham}(M,N)$. For each $k,r,s$ we have that
$$
\{\alpha_k,\alpha_r\}_{P_s} =0\;.
$$
Moreover, let $I_k={\rm Tr}N^k/k$. It satisfies
\begin{equation}\label{canonical_hierarchy}
N^*d I_k = d I_{k+1}\;\; .
\end{equation}
\end{proposition}
{\it Proof}. [after \cite{MM}] By using Proposition (\ref{consequences_pn_structure}) ($iii$) we get that 
$$
\{\alpha_k,\alpha_r\}_{P_s} = N^*{}^{k+r-2}\{\alpha,\alpha\}_{P_s} = 0 \;.
$$
The condition of vanishing torsion (\ref{Njienhuis_torsion}) is equivalent to state that for each vector field $v$ we have that $L_{N(v)}(N) = N L_v(N)$. So we compute for each 
vector field $v$
$$
\langle v, dI_{k+1}\rangle = \Tr(N^{k} L_v(N)) = \Tr(N^{k-1} L_{N(v)}(N)) = \langle N(v),dI_k\rangle = \langle v, N^*dI_k \rangle . \qed  
$$

\smallskip
We call the collection of hamiltonian $I_k$ the {\it canonical hierarchy}; they are in involution with respect to
all the Poisson structures $P_j$ of the hierarchy defined in $(i)$ of Proposition \ref{consequences_pn_structure}. We want to investigate when they are independent
and define a completely integrable model.

We will be interested in the case when $P=\omega^{-1}$ is the inverse of a symplectic form $\omega$. We call this case a {\it symplectic PN structure}. In this case PN structures are completely characterized by compatible
Poisson structures. Indeed, we have the following converse of Proposition \ref{consequences_pn_structure} ($i$).

\smallskip
\begin{lemma}
\label{compatible_poisson_structures_symplectic_case}
Let $\omega$ be a symplectic form and $\pi$ a Poisson tensor on $M$. If they are compatible, {\it i.e.} $[\pi,\omega^{-1}]=0$ then $(M,\omega^{-1},N=\pi\circ\omega)$
is a PN manifold.
\end{lemma}

{\it Proof}. Corollary 1.4 in \cite{vaisman}. \qed

\medskip


The algebroid morphism $P_1=\omega^{-1}:T^*_{P_j}M\rightarrow T_{N^{j-1}}M$ is invertible; in particular $T^*_{P_2}M=T^*_\pi M$ and $T_NM$ are isomorphic and the algebroid cohomology
$H(T_NM)$ coincides with the Lichnerowicz-Poisson cohomology of $P_2=\pi$. 

Let us consider now the eigenvalue problem for $N^*$. Since $N^*-\lambda=\omega\circ(\pi-\lambda\omega^{-1})$, with $\pi-\lambda\omega^{-1}$ being an antisymmetric operator, 
the eigenvalues of $N^*$ are doubly degenerate and the number of distinct eigenvalues can be at most $\dim M/2$.
The eigenspace corresponding to an eigenvalue $\lambda_x$ at $x\in M$ coincides with the kernel
of the antisymmetric operator $\pi-\lambda_x\omega^{-1}: T_x^*M\rightarrow T_xM$. 
We say that {\it $N$ is of maximal rank} if there exist an open dense $M_0$ where there are defined $\dim M/2$ independent functions $\lambda_i\in C^1(M_0)$ such that
$\lambda_i(x)$ is an eigenvalue of $N$ at $x\in M_0$. We call such functions $\lambda_i$ the {\it Nijenhuis eigenvalues}. The nondegeneracy condition is clearly
a statement about the complete integrability of the canonical hierarchy $\{I_k\}$ defined in Proposition \ref{PN_hierarchies}. Indeed, we easily see that 
$$
dI_1\ldots dI_n = \det(B) d\lambda_1\ldots d\lambda_n\;,
$$
where $n=\dim M/2$ and $B_{ik}=\frac{\partial I_k}{\partial\lambda_i}$ is the Vandermonde matrix, whose determinant is $\det B = \Pi_{i<j}(\lambda_j-\lambda_i)$.
Without loss of generality we can assume that on 
the open dense $M_0$ all eigenvalues are distinct and different from zero so that, if $A_{ij}=\lambda_j^i$, then also $\det A\not=0$ when evaluated on $M_0$. Indeed, one 
computes $\det A = (\Pi_i\lambda_i)  \Pi_{i<j}(\lambda_j-\lambda_i)$.

\begin{rem}{\rm
There is not an obvious way of defining such eigenvalue functions in full generality. The dense
$M_0$ can be disconnected and the ordering of the eigenvalues can introduce inessential singularities: the examples that we will discuss in the 
rest of the paper will clarify this issue. Moreover, there is no guarantee that they extend to the whole $M$; in our examples they will extend as continuous
functions.}
\end{rem}

\begin{proposition}\label{nijenhuis_eigenvalues_non_degenerate}
If $N$ is of maximal rank, the Nijenhuis eigenvalues satisfy the following equation
\begin{equation}\label{fundamental_eigenvalue_equation}
 N^*d\lambda_i = \lambda_i d\lambda_i \;.
\end{equation}
Moreover, the Nijenhuis eigenvalues are in involution with respect to all the Poisson structures $P_j$ and the Lie algebra structures $\g_j(M,N)$ defined on $\Omega^1_{ham}(M,N)$ are abelian for all $j$.
\end{proposition}
{\it Proof}. We compute from
(\ref{canonical_hierarchy}) that
$$
0 = N^*dI_k - d I_{k+1} = 2\sum_{i=1}^n\lambda^{k-1}_i (N^*d\lambda_i - \lambda_id\lambda_i)\;;
$$
relation (\ref{fundamental_eigenvalue_equation}) then follows because the nondegeneracy hypothesis implies that $\det A\not=0$ when evaluated on $M_0$ so that this linear 
system for $N^*d\lambda_i-\lambda_id\lambda_i$ admits only the zero solution. The eigenvalues are in involution because $I_k$ are in involution, as stated in Proposition 
\ref{PN_hierarchies}.

Moreover, since $d\lambda_i\not=0$, then $d\lambda_i$ is an eigenvector of $N^*$; if we denote with $V_{\lambda_i}\subset TM_0$ the eigenspace 
of $N$ corresponding to 
$\lambda_i$, then $\dim V_{\lambda_i}=2$ and
$$
TM_0 = \oplus_{i=1} V_{\lambda_i}\;.
$$
Let $v,w$ be vector fields parallel to $V_{\lambda_i}$; then we compute
\begin{eqnarray*}
0&=&T(N)(v,w)= [Nv,Nw]- N[Nv,w] - N[v,Nw] + N^2[v,w] \cr
&=& \lambda_i^2[v,w] - 2 \lambda_i N [v,w] + N^2[v,w] = (N-\lambda_i)^2[v,w]\;,
\end{eqnarray*}
so that, by the hypothesis of nondegeneracy, $[v,w]\in V_{\lambda_i}$, {\it i.e.} $V_{\lambda_i}$ is an involutive distribution. Let $\alpha\in\Omega^1_{ham}(M,N)$ and 
let $\alpha=\sum_i\alpha_i$ with $\alpha_i\in V_{\lambda_i}^*$. Moreover, since $V_{\lambda_i}$ is an involutive distribution then $d\alpha=dN^*\alpha=0$ implies 
$d\alpha_i=d(\lambda_i\alpha_i)=0$ so that
$\alpha_i=f_i(\lambda_i)d\lambda_i$. As a consequence $\{\alpha,\beta\}_{P_j}=0$ for all $\alpha,\beta\in\Omega^1_{ham}(M,N)$ and each $j$.
\qed

\bigskip
\bigskip

\section{Multiplicative integrability of the modular function}

We discuss in this section how a maximal rank PN structure defines a multiplicative integrable model. Let us consider two 
compatible Poisson structures, $\pi$ and $\omega^{-1}$, where $\omega$ is a symplectic form. We know from Lemma 
\ref{compatible_poisson_structures_symplectic_case} that the tensor $N=\pi\circ\omega$ is Nijenhuis and $(M,\omega^{-1},N)$ is a symplectic PN manifold. 
We will be interested in the Poisson geometry of the pencil $\pi_t=\pi+t \omega^{-1}$. By a straightforward
computation, it is easy to see that $N_t=\pi_t\circ\omega=N+t$ is a Nijenhuis tensor and that $(M,\omega^{-1},N_t)$ is a symplectic PN manifold for each $t\in\R$.
Since $\iota_{N_t^*}=\iota_{N^*} + t\deg$, where $\iota_{N^*}$ is the degree zero derivation defined before 
(\ref{Nijenhuis_algebroid_differential}) and $\deg(\nu)= k \nu$ for $\nu\in\Omega^k(M)$, then the algebroid differential is computed as 
$$d_{N_t}= [\iota_{N_t^*},d]= d_N + t d\quad\quad . $$
It is then clear that the space of hamiltonian forms does not depend on $t$. Given a hamiltonian form $\alpha\in\Omega^1_{ham}(M,N)$, the associated hierarchy
$\alpha_k(t)=(N^*+t)^k\alpha$ clearly depends on $t$. We assume in this section that $M$ is simply connected.

Let us first discuss the cohomological information of the hierarchy defined by a hamiltonian form. We know that the algebroid cohomology $H(T_{N_t}M)$ is isomorphic
to the Lichnerowicz-Poisson cohomology of $\pi_t$ by means of the invertible algebroid morphism $\omega:T_{N_t}M\rightarrow T^*_{\pi_t}M$.
A hierarchy of hamiltonian
forms is given by a hierarchy of functions $f_k\in C^\infty(M)$ satisfying
\begin{equation}\label{Lenart_hierarchy}
\omega^{-1}d f_{k+1} = \pi_t d f_k \;\quad\quad k\geq 1 \;. 
\end{equation}
We call a collection $\{f_k\}$ satisfying (\ref{Lenart_hierarchy}) a {\it Lenart hierarchy}.
As a consequence, $\sigma_{f_k}= \omega^{-1}df_k$ is a Poisson vector field for $\pi_t$ whose cohomology class 
is trivial if its hierarchy can be extended to one lower degree, this fact depending on $t$. 

The canonical hierarchy $I_k(t) = \Tr N^k_t/k$ extends down till $k=1$; it is shown in \cite{DaFe} that 
$\sigma_{I_1}$, independent on $t$, is the modular
vector field of $\pi_t$ with respect to the symplectic volume form; moreover, for those $t$ such that $\det N_t$ never vanishes, $N^*_tdI_0=dI_1$ where $I_0=\log\det N_t$.

It is clear from (\ref{fundamental_eigenvalue_equation}) that a Nijenhuis eigenvalue $\lambda$, if it 
extends to a global smooth function on $M$, defines a Lenart hierarchy. So each linear combination of the eigenvalues that extends to a smooth global function defines
a class in $H^1(T_{\pi_t}^*M)$: for instance $I_1 = \Tr N = \sum_i\lambda_i$.
 
Let $\G\equiv\G(M,\pi_t)$ be the ssc symplectic groupoid 
integrating $(M,\pi_t)$ and let $l,r$ denote the source and target maps. Let us denote with $h_\alpha\in C^\infty(\G)$ the groupoid one-cocycle integrating the Poisson vector
field $\omega^{-1}(\alpha)$ associated to $\alpha\in\Omega^1_{ham}(M,N)$ (see Appendix \ref{appendix_integrating_poisson_vector_fields} for the background). 
Let us consider a Lenart hierarchy $\{f_k\}$ and let $\{h_{df_k}\}$ be the corresponding cocycles. Since hamiltonian vector fields are integrated to trivial groupoid cocycles,
as a consequence of (\ref{Lenart_hierarchy}), we have that
\begin{equation}
\label{lenart_hierarchy_cocycles}
h_{df_{k+1}} = \partial^*(f_k) = l^*(f_k)- r^*(f_k)\quad ,
\end{equation}
where $\partial^*$ denotes the simplicial groupoid coboundary operator. Let us consider $\H(N_t)\subset\G(M,\pi_t)$ defined as
\begin{equation}
\label{PN_hamiltonian_subgroupoid}
\H(N_t) = \{\gamma\in\G(M,\pi_t)\ |\ h_{df}(\gamma)=0\ \forall df\in\Omega^1_{ham}(M,N_t)\}= \bigcap_{df\in\Omega^1_{ham}} \Ker h_{df}\;.
\end{equation}

Being defined as the intersection of the kernel of certain groupoid cocycles, then $\H(N_t)$ is a normal subgroupoid of $\G(M,\pi_t)$ so that we can define the 
quotient groupoid $\G(M,\pi_t)/\H(N_t)$. Let us define
\begin{equation}
\label{multiplicative_poisson_subalgebra}
\Omega^1_{ham}(\G(M,\pi_t))\equiv\{h_{df},l^*(f) |\ df\in\Omega^1_{ham}(M,N_t)\}\subset C^\infty(\G(M,\pi_t))\quad\quad .
\end{equation}

\smallskip
We prove the following result.

\begin{proposition}\label{pn_quotient_groupoid}
Let $M$ be a simply connected manifold equipped with two compatible Poisson structures $\omega^{-1}$ and $\pi$ and, for each $t\in\R$, let $\G(M,\pi_t)$ be
the ssc symplectic groupoid integrating $\pi_t=\pi+t\omega^{-1}$ (assumed to be integrable).
\begin{itemize}
\item[$i$)] $\Omega^1_{ham}(\G(M,\pi_t))$ is a Lie subalgebra of $C^\infty(\G(M,\pi_t))$ isomorphic to a central
$\R$-extension of
$$
\g_1(M,N_t)\ltimes\g_2(M,N_t)
$$
where $\g_i(M,N_t)$ are the Lie algebras introduced at the end of Proposition \ref{consequences_pn_structure} and $\g_1(M,N_t)$ acts on $\g_2(M,N_t)$ with 
$df\cdot dg = d\{f,g\}_{\omega^{-1}}$.
\item[$ii$)] For each $df\in\Omega^1_{ham}(M,N)$, $l^*(f)$ and $h_{df}$ descend to $\G(M,\pi_t)/\H(N_t)$. 
\item[$iii$)] The contour level set of $\Omega^1_{ham}(\G(M,\pi_t))$ inherits a topological groupoid structure, that we denote as $\L(N_t)$, such that
\begin{equation}\label{groupoid_morphisms}
\begin{tikzcd}
\G(M,\pi_t) \arrow{r} \arrow{d}
&\L(N_t) \\
\G(M,\pi_t)/\H(N_t) \arrow{ru} & 
\end{tikzcd}
\end{equation}
is a commutative diagram of surjective groupoid morphisms.
\end{itemize}
\end{proposition}
{\it Proof}. $i$) 
If we apply (\ref{bracket_groupoi_cocycles}) to $\sigma_{g_i}=\omega^{-1}(dg_i)$, with $dg_i\in\Omega^{1}_{ham}(M,N)$, we get 
$$
\{h_{dg_1},h_{dg_2}\} = h_{d\{g_1,g_2\}_{\omega^{-1}}}
$$
with $d\{g_1,g_2\}_{\omega^{-1}}\in\Omega^1_{ham}(M,N)$, since hamiltonian forms are closed with respect to both brackets defined by $\pi_t$ and $\omega^{-1}$, as 
stated in Proposition \ref{nijenhuis_eigenvalues_non_degenerate}. Analogously let us apply (\ref{bracket_cocycle_pullback}) to $\omega^{-1}dg$ and to $l^*(f)$, with
$df,dg\in\Omega^{1}_{ham}$ and get
$$
\{h_{dg},l^*(f)\} = l^*(\omega^{-1}dg(f)) = l^*{\{g,f\}_{\omega^{-1}}}.
$$
Finally, if $df,dg\in\Omega^1_{ham}(M,N)$ then
$$
\{l^*(f),l^*(g)\}= l^*\{f,g\}_{\pi_t}
$$
since $l^*$ is a Poisson morphism.

$ii$) Indeed, let $\gamma'= \xi \gamma \xi'$, with $\xi,\xi'\in \H(N_t)$. Then we compute
$$
h_{df}(\gamma')= h_{df}(\xi)+h_{df}(\gamma) + h_{df}(\xi')= h_{df}(\gamma)
$$
and, if $N^*df = df_1$,
$$
l^*(f)(\gamma') = f(l(\gamma')) = f(l(\xi)) = f(r(\xi))+h_{df_1}(\xi) = f(l(\gamma)) = l^*(f)(\gamma)\;,
$$
where in the second equality we used (\ref{lenart_hierarchy_cocycles}).

$iii$) Let $\gamma,\gamma'\in\G(M,\pi_t)$ such that $h_{df}(\gamma)=h_{df}(\gamma')=h_{df}$ and $f(l(\gamma))=f(l(\gamma'))=f$ for each Lenart hierarchy $f=\{f_k\}$. From 
(\ref{lenart_hierarchy_cocycles}) we see that
$$
f_k(r(\gamma)) = f_k(l(\gamma)) - h_{df_{k+1}}(\gamma) = f_k-h_{df_{k+1}}= f_k(r(\gamma'))\;.
$$
Moreover, if $\gamma,\gamma'\in\G_2(M,\pi_t)$ then $h_{df_k}(\gamma\gamma')=h_{df_k}(\gamma)+h_{df_k}(\gamma')$ and $l^*(f_k)(\gamma\gamma')=l^*(f_k)(\gamma)$, {\it i.e.} they
are indepedent on the choice of $\gamma,\gamma'$ on the contour level set. \qed

\medskip

\begin{rem}{\rm
The kernels of the two surjective groupoid morphisms appearing in (\ref{groupoid_morphisms}) coincide with $\H(N_t)$. This doesn't imply that the diagonal morphism
is an isomorphism, because in groupoids, differently than in groups, surjective morphisms are characterized by their kernels only if they are base and piecewise surjective (see
\cite{MK}).}
\end{rem}

\smallskip
\begin{example}
{\rm Let us consider the completely degenerate case. Let $\omega$ be a symplectic form and let us choose $N={\rm id}$ so that all the Poisson structures of the PN hierarchy coincide
with $\omega^{-1}$. It is clear that hamiltonian forms coincide with exact forms; moreover, both $\L(N)$ and $\G(M,\omega^{-1})/\H(N)$ coincide with $\G(M,\omega^{-1})$ itself.}
\end{example}

\smallskip
\begin{corollary}
\label{multiplicative_integrable_model}
If $N_t$ is of maximal rank
then $\Omega^1_{ham}(\G(M,\pi_t))$ is abelian and defines a multiplicative integrable model.
\end{corollary}
{\it Proof}. It is abelian since both $\g_1(M,N)$ and $\g_2(M,N)$ are abelian Lie algebras, as seen in Proposition \ref{nijenhuis_eigenvalues_non_degenerate}. 
On the open dense $\tilde{M}_0=\{m\in M_0|\, \lambda_t(m)\not=0\}$ $\pi_t$ is of maximal rank and 
$\G(M,\pi_t)|_{\tilde{M}_0}=\tilde{M}_0\times\tilde{M}_0$; when restricted to this open dense, the hamiltonians in involution are just the tensor product of the hamiltonian
on $\tilde{M}_0$ and so they are independent. \qed

\smallskip
In \cite{BCQT} we called such a system a {\it multiplicative integrable model}. Recall that the canonical Lenart hierarchy $I_k=\Tr(N^k_t)/k$ describes the modular class of $\pi_t$;
the hamiltonian form $dI_1 = d\Tr(N_t)$ lifts to $h_{dI_1}$, 
that we call the {\it modular function}. We can rephrase the above result by saying that the modular function $h_{dI_1}$ is multiplicatively integrable. 

\smallskip
Let us assume now that the Nijenhuis tensor is of maximal rank and that the eigenvalues exist as global continuous functions while they are smooth only on the dense open 
$\widetilde{M}_0$. Then the space of identities of  
${\cal L}(N_t)$ is the image of the Nijenhuis eigenvalues that
we denote as $\Delta(N_t)=\Delta(N)+t\subset\R^n$, where $\dim M=2n$. By means of an obvious redefinition, we consider that the space of identities of $\L(N_t)$ is $\Delta(N)$ for all $t$, that we 
call the {\it bihamiltonian polytope}. The image of $\widetilde{M}_0$ inside $\Delta(N)$ is obviously indepedent on $t$ and is denoted as $\Delta(N)_0$.
On $\G(M,\pi_t)|_{\widetilde{M}_0}$ we can lift the Poisson vector field $\sigma_{\lambda_i}$ associated with the eigenvalue $\lambda_i+t$ to the cocycle 
$h_{d\lambda_i}$;
$\{\lambda,h_{d\lambda}\}$ is a chart of continuous coordinates for $\L(N_t)|_{\Delta(N)_0}$. The groupoid 
structure maps in terms of these coordinates become particularly simple.
Indeed, if we restrict to the locus $M_0\subset\widetilde{M_0}$, depending on $t$, where the eigenvalues $\lambda_i(t)=\lambda_i+t\not=0$, we see that
$$
\pi_td\log(\lambda_i+t)= \pi_td\lambda_i/(\lambda_i+t) = \omega^{-1}(N^*+t)d\lambda_i/(\lambda_i+t) = \omega^{-1}d\lambda_i = \sigma_{\lambda_i}\;,
$$
{\it i.e.} $\log(\lambda_i+t)$ is a local hamiltonian for $\sigma_{\lambda_i}$ so that 
\begin{equation}\label{lift_eigenvalues}
h_{d\lambda_i} = \partial^* \log(\lambda_i+t) \;.
\end{equation}
If $\gamma\in\G(M,\pi_t)|_{M_0}$ is such that $l^*(\lambda)(\gamma) = \lambda$ and 
$h_{d\lambda}(\gamma) = h_{\lambda}$ then we get that
\begin{equation}\label{groupoid_map_from_eigenvalue}
l(\lambda,h_{d\lambda})=\lambda,\quad r(\lambda,h_{d\lambda})= -t + e^{h_{d\lambda}}(\lambda+t),\quad (\lambda,h_{d\lambda})(\lambda',h_{d\lambda}')=(\lambda,h_{d\lambda}+
h_{d\lambda}')\;.
\end{equation}
These formulas extend wherever $h_{d\lambda_i}$ is defined, at least, as a continuous function. Remark that the lift of the Poisson vector field $\omega^{-1}d\lambda_i$ 
defined in (\ref{lift_eigenvalues}) extends whenever $\lambda_i+t\not= 0$, while the lift described in (\ref{lift_poisson_vector_cotangent_path}) is defined whenever
$d\lambda_i\not =0$. It is in general relavant to understand the maximal extension of $h_{\lambda_i}$, in particular if it extends at least as a global continous cocycle.
This is the case for instance if $|t|$ is big enough such that $\pi_t$ is non degenerate and 
$\lambda_i+t\not= 0$ everywhere. In this case we can use $l(\lambda,h_{d\lambda})$ and $r(\lambda,h_{d\lambda})$ as coordinates and $\L(N_t)$ is the pair groupoid 
$\Delta(N)\times\Delta(N)$.

Let us consider the action of the additive group 
$\R^n$ on $\R^n$ defined by $r$ above: $h\in\R^n$ acts on $\lambda\in\R^n$ as
\begin{equation}\label{abelian_action}
h(\lambda) = r(h,\lambda) = -t + e^{h}(\lambda + t)\:. 
\end{equation}
Then the restricted groupoid $\L(N_t)|_{\Delta(N)_0}$ can be described as a subgroupoid of the action groupoid restricted to $\Delta(N)_0\subset\R^n$
$$
\L(N_t)|_{\Delta(N)_0} \subset \R^n\rtimes\R^n|_{\Delta(N)_0}\;.
$$

\smallskip
The global description of $\L(N_t)$ is an important point. Moreover, it is crucial to clarify if
the projection $\G(M,\pi_t)\rightarrow {\cal L}(N_t)$ has lagrangian fibres so that the groupoid of Bohr-Sommerfeld 
leaves can be defined according to the general framework proposed in \cite{BCQT}.

\section{PN structures of maximal rank on $Gr(k,n)$}

We give in this section an example of multiplicative integrable model based on the construction of the maximal rank PN structures 
on compact hermitian symmetric spaces discussed in \cite{BQT}. 
Let $M$ be a compact hermitian G-symmetric space. The compatible Poisson structures are $(P,\pi)$, where $P=\omega_{kks}^{-1}$ is the inverse of the Kirillov-Kostant-Souriau 
symplectic form, defined when $M$ is seen as a coadjoint orbit, and $\pi$ 
is the so-called {\it Bruhat-Poisson structure}, obtained as quotient of the standard Poisson-Lie group structure on $G$. Their compatibility
$$
[\pi,\omega_{kks}^{-1}] = 0
$$
has been observed first in \cite{KRR}. The associated pencil $\pi_t = \pi +t \omega_{kks}^{-1}$ defines a family of
Poisson homogeneous spaces of the Poisson-Lie group $G$. 

This construction will be sketched here only for $G=U(n)$ and $M=Gr(k,n)$.
Let us consider $Gr(k,n)$ as the adjoint $U(n)$-orbit through 
$$\rho={\rm diag}(\underbrace{i,\ldots,i}_{k},\underbrace{0,\ldots,0}_{n-k})\in\u(n)\;.$$
Here we clearly identify $\u(n)$ with $\u(n)^*$ via the trace and this fixes the Kirillov-Kostant-Souriau symplectic form $\omega_{kks}$ that is invariant under the $U(n)$ transformations. The Iwasawa decomposition of $\sl(n,\C)$ defines a Manin triple and so
a bialgebra structure on $\su(n)$ trivially extended to $\u(n)$, that we denote with $\delta_{\u(n)}:\u(n)\rightarrow\Lambda^2\u(n)$. It integrates to a Poisson-Lie group 
structure on $U(n)$, that we call the {\it standard Poisson-Lie group structure}. We skip explicit details that can be found in \cite{BQT}.
The Poisson tensor $\pi_{U(n)}$ is projectable with respect to the quotient map $U(n)\rightarrow Gr(k,n)$. We denote with $\pi$ the induced Poisson structure and 
we call it {\it Bruhat-Poisson structure}. By construction $\pi$ is just covariant, {\it i.e.}
$L_X(\pi)=\sigma(\delta_{\u(n)}(X))$, for each $X\in\u(n)$ and $\sigma:\u(n)\rightarrow {\rm Vect}(Gr(k,n))$ denotes the infinitesimal action.
It was shown in \cite{KRR} that $\pi$ and $\omega_{kks}^{-1}$ are compatible so that by Lemma \ref{compatible_poisson_structures_symplectic_case}
they define a PN structure on $Gr(k,n)$ with Nijenhuis tensor $N=\pi\circ\omega_{kks}$.

In \cite{BQT} it has been shown that this symplectic PN structure is of maximal rank. This means that there exist $\dim Gr(k,n)/2=k(n-k)$ Nijenhuis eigenvalues
that are independent on an open dense subset. 
%
Let us consider the following chain of subalgebras
\bea
\u(n)\supset \u(n-1)\ldots \supset \u(1)\label{chain_subalgebra_su}
\eea
acting in a hamiltonian way on $Gr(k,n)$ with moment map $\mu_{\u(s)}$ with $s=1,\ldots,n$. Since we identify $\u(s)$ with $\u(s)^*$, $\mu_{u(s)}$ can be identified as a $s\times s$ matrix;
in particular, if we choose to embed $\u(s)$ in the upper left minor of $\u(n)$, then $\mu_{\u(s)}$ is just this minor of $\mu_{\u(n)}$. It is proven in 
\cite{BQT} that each non constant eigenvalue of $\mu_{\u(s)}$ defines a Nijenhuis eigenvalue of $N$. These variables 
are the well known {\it Gelfand-Tsetlin} variables
that define a completely integrable model on all partial flag manifolds. In particular they are independent and this proves that the PN structure is of maximal rank.

Let $\G(Gr(k,n),\pi_t)$ be the symplectic groupoid integrating the Poisson structure $\pi_t = \pi + t\omega_{kks}^{-1}$. From the previous section we know how to 
define a multiplicative integrable model. Let $\L(N_t)$ be the topological groupoid defined in Proposition \ref{pn_quotient_groupoid} ($iii$).
Since the Gelfand-Tsetlin variables are globally continuous functions, the space of identities of $\L(N_t)$ is the so-called Gelfand-Tsetlin polytope 
${\cal C}_{GT}(k,n)\subset\R^{k(n-k)}$, defined as the space of independent solutions of the GT inequalities
$$
0\leq \lambda^{(s)}_i\leq\lambda^{(s+1)}_i\leq\lambda^{(s)}_{i+1}\leq 2 \quad\quad i=1,\ldots n-s,\; s=1,\ldots n-1\;.
$$
If $t\not\in[-2,0]$ then $\lambda_i+t\not=0$ everywhere and $\pi_t$ is non degenerate. In this case $\L(N_t)$ is just the pair groupoid 
${\cal C}_{GT}(k,n)\times{\cal C}_{GT}(k,n)$, as it has been discussed at the end of the previous section.

When $t\in[-2,0]$, we don't have a global description of $\L(N_t)$. In particular,
since the GT variables, when restricted to the boundary of the polytope, fail to be smooth, we don't know if they lift to global continuous groupoid cocycles of
$\G(Gr(k,n),\pi_t)$. 

The global description of $\L(N_t)$ has been done so far only in the case $k=1$, {\it i.e.} for the complex projective space $\C P_{n}=Gr(1,n+1)$. In this case the 
Gelfand-Tsetlin variables are global smooth functions. 
Indeed, for each $s\leq n$ the unique non constant eigenvalue of the moment map $\mu_{\u(s)}$ of the 
chain (\ref{chain_subalgebra_su}) is the hamiltonian of the vector field of the action of
$$
H_s=2i\ {\rm diag}(\underbrace{1,\ldots,1}_{s},\underbrace{0,\ldots,0}_{n+1-s})\in\u(n+1)\,.
$$
With the normalization of the Poisson tensors specified in \cite{BQT}, $\pi_t$ is non degenerate for $t\not\in[-2,0]$; 
in this case then the first Poisson cohomology group vanishes. If $t\in[-2,0]$, the Poisson cohomology class of the vector fields of the Cartan 
action is non trivial and we can identify the image of $\g_1(\C P_{n},\pi_t)$ in the cohomology of $\pi_t$ with $\t_n$, the Cartan subalgebra of $\su(n+1)$. 

The Gelfand-Tsetlin polytope is in this case just the simplex $\Delta_{n}=\{\lambda\in\R^{n}| 0\leq\lambda_1\leq\lambda_2\leq\ldots\lambda_{n}\leq 2\}$. 
The eigenvalues and their lifted cocycles are then a set of global coordinates on $\L(N_t)$ with the structure maps defined in 
(\ref{groupoid_map_from_eigenvalue}). For $t\not\in[-2,0]$, $\L(N_t)$ is the pair
groupoid $\Delta_n\times\Delta_n$. For $t\in[-2,0]$, the following description is the content of Proposition 6.1 of \cite{BCQT}. Recall the action of $\R^n$ on $\R^n$
defined in (\ref{abelian_action}) and let $\R^n\rtimes\R^n$ denote the action groupoid.
 
\smallskip
\begin{proposition}\label{multiplicative_integrable_model_cpn}
For $t\in (-2,0)$, $\L(N_t)$ is isomorphic to
the following wide subgroupoid of $\R^n\rtimes\R^n|_{\Delta_n}$,
$$ \{(\lambda,h_\lambda)\in\R^n\rtimes\R^n|_{\Delta_n}~|~ \lambda_i=\lambda_{i+1}=-t \implies~ h_{\lambda_i}=h_{\lambda_{i+1}}\}~.$$
For $t\in \{-2,0\}$, $\L(N_t)$ is isomorphic to the wide subgroupoid
$$\{(\lambda,h_\lambda)\in\R^n\rtimes\R^n|_{\Delta_n}~|~ \lambda_i=-t \implies~ h_{\lambda_i}=0\}~.$$
\end{proposition}

In \cite{BCQT} it has been shown that the quotient map $\G(\C P_{n-1},\pi_t)\rightarrow \L(N_t)$ described in Proposition
\ref{pn_quotient_groupoid} has lagrangian fibres and the Bohr-Sommerfeld groupoid has been computed. The resulting groupoid admits a unique Haar system so that
the quantization program can be completed till the definition of the convolution algebra. We refer to \cite{BCQT} for all details.

The general case of Grassmannians and the other compact hermitian symmetric spaces will be analyzed in another publication. We can expect that a major difference will appear due 
to the non toric nature of these manifolds. In particular, the Nijenhuis eigenvalues are just global continuous functions and the description of the groupoid 
$\L(N_t)$ given in (\ref{groupoid_map_from_eigenvalue}) is valid only when restricted on (the image of) the open dense subset where the eigenvalues are smooth. 
The main question to understand is if the groupoid quotient map has lagrangian fibres and allows the definition of Bohr-Sommerfeld leaves.

\bigskip
\bigskip

\appendix

\section{From Poisson vector fields to groupoid cocycles}\label{appendix_integrating_poisson_vector_fields}
We recall in this Appendix basic facts about the lift of Poisson vector fields to groupoid cocycles.

Let $P$ be an integrable Poisson structure on $M$ and let $\G=\G(M,P)$ be the source simply connected (ssc) Lie groupoid integrating it. 
We denote with $l_\G$ and $r_\G$ the source and target maps and with $\G_s$ the space of strings of $s$-composable elements of $\G$, where $\G_0=M$ and $\G_1=\G$.
The face maps are $d_i:\G_s\rightarrow \G_{s-1}$, $i=0,\ldots s$, defined for $s>1$ as
\begin{equation}\label{face_maps}
d_i(\gamma_1,\ldots\gamma_s) = \left\{ \begin{array}{ll} (\gamma_2,\ldots\gamma_s)& i=0\cr
(\gamma_1,\ldots\gamma_i\gamma_{i+1}\ldots) & 0<i<s\cr
(\gamma_1,\ldots\gamma_{s-1})& i=s\end{array}\right.
\end{equation}
and for $s=1$ as $d_0(\gamma)=l_\G(\gamma)$, $d_1(\gamma)=r_\G(\gamma)$. The simplicial coboundary operator
$\partial^*:\Omega^k(\G_{s})\rightarrow\Omega^k(\G_{s+1})$ is defined as
$$
\partial^*(\omega) = \sum_{i=0}^s (-)^i d_i^*(\omega)  \;,
$$
and ${\partial^*}^2=0$.
The cohomology of this complex for $k=0$ is the real valued groupoid cohomology; $s$-cocycles are denoted as $Z^s(\G,\R)$.

Let $\sigma\in{\rm Vect}^1(M)$ be a Poisson vector field, {\it i.e.} it is closed under the algebroid differential $d_P$. We want to show that it can
be lifted to a groupoid one-cocycle. It is very useful to use the construction of the symplectic groupoid $\G$ as symplectic reduction from $T^*M^{[0,1]}$, where 
$M^{[0,1]}$ is the path space, 
as proven in \cite{CatFel}. Indeed, if $(X,\eta)\in T^*M^{[0,1]}$, where $X\in M^{[0,1]}$ and $\eta\in\Gamma(X^*T^*M)$ then $\G$ can be constructed as 
the symplectic reduction with respect to the constraint
$$
\dot{X} + P(\eta) = 0\;.
$$
The source and target map are then $l_\G(X,\eta)=X(0)$ and $r_\G(X,\eta)=X(1)$ and multiplication is by concatenation of paths. We can lift $\sigma$ to a 
groupoid one-cocycle
\begin{equation}\label{lift_poisson_vector_cotangent_path}
h_{\sigma}(X,\eta)= \int_0^1 \langle \sigma,\eta\rangle dt\;.
\end{equation}
Indeed, it can be checked that since $d_P(\sigma)=0$, the value of $h_s(X,\eta)$ doesn't change if deform $(X,\eta)$ by a cotangent homotopy so that $h_\sigma\in C^\infty(\G)$. 
Since groupoid multiplication is just concatenation of paths, then it is clear that $h_{\sigma}$ is a groupoid one-cocycle, 
$\partial^*(h_\sigma)=0$. Moreover, it can be checked that if $\sigma=d_P(f)$ for some 
$f\in C^\infty(M)$ then
$$
h_\sigma = \partial^*(f)\;.
$$
Let us consider $\sigma_i$ be Poisson vector fields and let $h_{\sigma_i}$ the lifted 
groupoid cocycles. The symplectic structure on $\G$ comes from the symplectic reduction of the canonical symplectic form on $T^*M^{[0,1]}$ so that 
we can compute 
\begin{equation}\label{bracket_groupoi_cocycles}
\{h_{\sigma_i},h_{\sigma_j}\} = \langle \Omega_{M^{[0,1]}}^{-1},d h_{\sigma_i}\wedge d h_{\sigma_j}\rangle = \int_0^1 
(\frac{\delta h_{\sigma_i}}{\delta X^\mu}\frac{\delta h_{\sigma_j}}{\delta \eta_\mu}-
\frac{\delta h_{\sigma_j}}{\delta X^\mu}\frac{\delta h_{\sigma_i}}{\delta \eta_\mu})dt= 
h_{[\sigma_i,\sigma_j]}\;.
\end{equation}
Analogously, if $f\in C^\infty(M)$ then $l^*_\G(f)(X,\eta)= f(X(0))$; we then compute
\begin{equation}\label{bracket_cocycle_pullback}
\{h_{\sigma_i},l^*_\G(f)\} = l^*_\G(\sigma_i(f))\;.
\end{equation}


\begin{thebibliography}{6666}

\bibitem{BCST2} Bonechi F., Ciccoli N., Staffolani N., Tarlini M., {\it The quantization of the symplectic groupoid of the standard Podle\`s sphere},
Journal of Geometry and Physics {\bf  62}, 1851--1865, (2012).

\bibitem{BCQT} Bonechi F., Ciccoli N., Qiu J., Tarlini M., {\it Quantization of Poisson Manifolds from the Integrability of the Modular Function},
Commun. Math. Phys. {\bf  331}, 851--885, (2014).

\bibitem{BQT} Bonechi F., Qiu J., Tarlini M., {\it Complete integrability from Poisson-Nijenhuis structures on compact hermitian symmetric spaces}, 
arXiv:1503.07339 [math.SG]

\bibitem{CatFel} A.S. Cattaneo and G. Felder, {\it Poisson sigma models and symplectic groupoids}.
In ``Quantization of Singular Symplectic Quotients'', (ed. N. P. Landsman, M. Pflaum,
M. Schlichenmeier), Progress in Mathematics, {\bf 198} (2001) 41--73.

\bibitem{DaFe} P. A. Damianou and R. L. Fernandes, \emph{Integrable hierarchies and the modular class,} Ann. Inst. Fourier {\bf 58}, (2008) 107--137.

\bibitem{Hawkins}
E. Hawkins, {\it A groupoid approach to quantization}. J. Symplectic Geom., {\bf 6} (2008)
 61--125  [arXiv:math.SG/0612363].
 
\bibitem{KRR}  S. Khoroshkin, A. Radul, V. Rubtsov, {\it A family of Poisson structures on hermitian symmetric spaces}. Commun. Math. Phys.
{\bf 152}, 2, (1993), 299-315. 

\bibitem{KS} Kosmann-Schwarzbach Y., {\it The Lie bialgebroid of a Poisson-Nijenhuis manifold}. Letters in Mathematical Physics (1996) {\bf 38}, 4, 421 -428 

\bibitem{MK} Mackenzie K., {\it Lie Groupoids and Lie algebroids in Differential Geometry}. London Mathematical Society, Lecture Note Series, 124, Cambridge University Press,
(1987).

\bibitem{MKS} F. Magri and Y. Kosmann-Schwarzbach, {\it Poisson-Nijenhuis structure}, Annales de l'I.H.P. Physique théorique, (1990), {\bf 53}, 1, 35-81.

\bibitem{MM} F. Magri and C. Morosi {\it A geometrical characterization of Hamiltonian systems through the theory of Poisson-Nijenhuis manifolds}.
Universit\`a di Milano-Bicocca, Quaderni di Matematica n.3/2008

\bibitem{vaisman} I. Vaisman, {\it The Poisson-Nijenhuis manifolds revisited}. Rend. Sem. Mat. Univ. Poi. Torino
Voi. 52, 4 (1994)


\end{thebibliography}
\end{document}